\documentstyle[12pt]{amsart}
\let\le\undefined
\let\ge\undefined
\newsymbol\le 1336    
\newsymbol\ge 133E    
\newcommand{\ds}{\dots}

\newcommand{\dsb}{\dotsb}
\newcommand{\cd}{\cdot}
\def\ot{\DOTSB\otimes}

\def\rarrow{\DOTSB\longrightarrow}
\def\larrow{\DOTSB\longleftarrow}
\newcommand{\maps}{\longmapsto}
\newcommand{\lrarrow}{\DOTSB\,\relbar\joinrel\relbar\joinrel\rightarrow\,}
\newcommand{\llarrow}{\DOTSB\,\leftarrow\joinrel\relbar\joinrel\relbar\,}
\newcommand{\trarrow}{\;\widetilde\longrightarrow\;}
\renewcommand{\:}{\colon}
\newsymbol\k 207C     
\newcommand{\ind}{}
\newenvironment{thm}[1]{\smallskip\ind{\sf #1.}\sl}{\smallskip}
\newenvironment{thmsec}[1]{\ind{\sf #1.}\sl}{\smallskip}
\newenvironment{rem}[1]{\smallskip\ind{\sf #1}}{\smallskip}
\newenvironment{remsec}[1]{\ind{\sf #1}}{\smallskip}
\newcommand{\pr}[1]{\ind{\it #1\/}}

\newcommand{\Section}[1]{\bigskip\section{#1}\medskip}

\newcommand{\coker}{\operatorname{coker}}
\newcommand{\id}{\operatorname{id}}
\newcommand{\chr}{\operatorname{char}}
\newcommand{\Hom}{\operatorname{Hom}}
\newcommand{\Tor}{\operatorname{Tor}}
\newcommand{\Ext}{\operatorname{Ext}}
\newcommand{\sub}{\subset}
\newcommand{\bop}{\bigoplus\nolimits}
\newcommand{\bcap}{\bigcap\nolimits}

\newcommand{\bs}{\{\,}
\newcommand{\es}{\,\}}
\newcommand{\lan}{\langle}
\newcommand{\ran}{\rangle}

\newcommand{\F}{{\Bbb F}}
\newcommand{\Z}{{\Bbb Z}}
\newcommand{\T}{{\Bbb T}}
\newcommand{\N}{{\Bbb N}}
\newcommand{\Q}{{\Bbb Q}}
\newcommand{\R}{{\Bbb R}}
\newcommand{\X}{{\cal X}}
\newcommand{\W}{{\cal W}}
\renewcommand{\a}{\alpha}
\renewcommand{\c}{\gamma}
\renewcommand{\o}{\omega}
\renewcommand{\O}{\Omega}
\renewcommand{\d}{\partial}
\newcommand{\e}{\varepsilon}
\newcommand{\D}{\Delta}
\newcommand{\Gal}{\operatorname{Gal}}
\newcommand{\oF}{\,\overline{\!F}}
\newcommand{\KM}{K^{\operatorname{M}}}
\newcommand{\q}{{\operatorname{q}}}
\newcommand{\gr}{{\operatorname{gr}}}

\begin{document}
\rightline{\scriptsize\hfill
Preprint alg-geom/9507010}
\vspace{0.75cm}
\title{Koszul Duality and Galois Cohomology}
\author{Leonid Positselski}
\address{Independent University of Moscow}
\email{posic@@ium.ips.ras.ru,\, posic@@math.harvard.edu}
\author{Alexander Vishik}
\address{Harvard University}
\email{vishik@@math.harvard.edu}
\maketitle

\section*{Introduction}
\smallskip

 Let $F$ be a field, $\oF$ be its (separable) algebraic closure, and
$G_F=\Gal(\oF/F)$ be the absolute Galois group.
 Let $l\ne\chr F$ be a prime number; assume that $F$ contains a $l$-root
of unity $\zeta$.
 In this case, the Kummer pairing
  $$
   \kappa\:G_F\times F^*\rarrow \F_l,
   \qquad \kappa(g,a)=s
   \text{ \ if \ }  g(b)=\zeta^sb
   \text{ \ for \ } b=\sqrt[\uproot2\leftroot1 l] a\in\oF
  $$
defines an isomorphism
$F^*/(F^*)^l\trarrow H^1(G_F,\F_l)$.
 The Milnor K-theory ring $\KM(F)$ is a skew-commutative quadratic algebra
over $\Z$ generated by $\KM_1(F)=F^*$ with the Steinberg relations
$\{a,1-a\}=0$.
 It is not difficult to show that the Kummer map can be extended
to an algebra homomorphism
  $$
   \KM(F)\ot\F_l\lrarrow H^*(G_F,\F_l),
  $$
which is known as the {\it Galois symbol}, or the
{\it norm residue homomorphism}.
 The well-known {\it Bloch--Kato conjecture\/} claimes that it is an
isomorphism.
 It was proved by A.~Merkurjev and A.~Suslin~\cite{MS1,MS2} and
M.~Rost~\cite{Ros} that this is true in degree~$2$ and
for~$l=2$ in degree~$3$.
 The aim of this note is to show that the whole conjecture follows
from its low-degree part provided the quadratic algebra $\KM(F)\ot\F_l$ is
{\it Koszul\/} (see section~2 for the definition).
 We will assume that $F$ has no algebraic extensions of degree
relatively prime to $l$.

\begin{thm}{Theorem}
 Let $H=H^*(G,\F_l)$ be the cohomology algebra of a pro-$l$-group $G$.
 Assume that
  \begin{enumerate}
    \item $H^2$ is generated by $H^1$;
    \item in the subalgebra generated by $H^1$ in $H$, there are no
          nontrivial relations of degree~3;
    \item the quadratic algebra defined by $H^1$ and $H^2$ is Koszul.
  \end{enumerate}
 Then the whole algebra $H$ is quadratic.
\end{thm}

 Actually, it is not essential that we consider the group cohomology here;
the theorem is valid for any pro-nilpotent algebra.
 We will use the language of coalgebras in this paper in order to avoid
dealing with projective limits and dualizations.

 We are grateful to V.~Voevodsky for stating the problem and numerous
stimulating discussions and to J.~Bernstein who pointed out to us
the necessity of using coalgebras in the Koszul duality.
 The first author is pleased to thank Harvard University for its
hospitality which made it possible for this work to appear.

\Section{Bar Construction}

\subsection{The cohomology of augmented coalgebras}
 A {\it coalgebra\/} is a vector space $C$ over a field $\k$ equipped
with a comultiplication map $\D\:C\rarrow C\ot C$ and a counit map
$\e\:C\rarrow\k$ satisfying the conventional associativity and counit
axioms.
 An {\it augmented coalgebra\/} is a coalgebra $C$ equipped with a
coalgebra homomorphism $\c\:\k\rarrow C$.
 The cohomology algebra of an augmented coalgebra $C$ is defined as
the $\Ext$-algebra $H^*(C)=\Ext^*_C(\k,\k)$ in the category of left
$C$-comodules,  where $\k$ is endowed with the comodule structure by means
of $\c$.
 We will calculate this cohomology using the following explicit
comodule resolution
  $$
   \k\lrarrow C\lrarrow C\ot C^+\lrarrow C\ot C^+\ot C^+\lrarrow\dsb,
  $$
where $C^+=\coker(\c)$, the differential is
  $$
   d(c_0\ot\ds\ot c_n)=\sum_{i=0}^n (-1)^{i-1} c_0\ot\ds\ot
   \D(c_i)\ot\ds\ot c_n,
  $$
and the coaction of $C$ is through the left components of these tensors.
 It is easy to check that the comodules $C\ot W$ are injective, the
differential is well-defined, and the operator
  $$
   h\: c_0\ot c_1\ot\ds\ot c_n \maps
    \e(c_0)\sigma(c_1)\ot c_2\ot\ds\ot c_n,
  $$
where $\sigma\:C^+\rarrow C$ is the splitting along $\e$,
provides a $\k$-linear contracting homotopy, thus this is a resolution.
 Applying the functor $\Hom_C(\k,\cd\,)$, we obtain
  $$
   H^*(C)=H^*(\,\k\rarrow C^+\rarrow C^+\ot C^+\rarrow\dsb\,),
  $$
where the differential is given by the same formula and the
multiplication on $H^*(C)$ is induced by the evident multiplication
  $$
   (c_1\ot\ds\ot c_i)\cd(c_{i+1}\ot\ds\ot c_{i+j})=
   c_1\ot\ds\ot c_i\ot c_{i+1}\ot\ds\ot c_{i+j}
  $$
on this cobar-complex.

\subsection{The homology of augmented algebras}
 An augmented algebra $A$ is an associative algebra over a field $\k$
endowed with an algebra homomophism $\a\:A\rarrow\k$.
 The homology coalgebra of an augmented algebra $A$ is by the definition
$H_*(A)=\Tor_A(\k,\k)$, where the left and right module structures on
$\k$ are defined by means of~$\a$.
 We will calculate it using the following explicit bar-resolution of the
left $A$-module~$\k$
  $$
   \k\llarrow A\llarrow A\ot A_+\llarrow A\ot A_+\ot A_+\llarrow\dsb,
  $$
where $A_+=\ker(\a)$ and
  $$
   \d(a_0\ot\ds\ot a_n)=\sum_{i=1}^n (-1)^i a_0\ot\ds\ot
   a_{i-1}a_i\ot\ds\ot a_n.
  $$
 It is easy to check that the operator
  $$
   h\: a_0\ot a_1\ot\ds\ot a_{n-1} \maps
   1\ot (a_0-\a(a_0))\ot a_1\ot\ds\ot a_{n-1}
  $$
provides a $\k$-linear contracting homotopy.
 Applying the functor $\k\ot_A\cd\,$, we obtain
  $$
   H_*(A)=H_*(\,\k\larrow A_+\larrow A_+\ot A_+\larrow\dsb\,)
  $$
and the coalgebra structure on $H_*(A)$ is induced by the evident
coalgebra structure
  $$
   \D(a_1\ot\ds\ot a_n)=\sum_{i=0}^n
   (a_1\ot\ds\ot a_i)\ot(a_{i+1}\ot\ds\ot a_n)
  $$
on this bar-complex.

\Section{Koszul Duality}

 By a graded algebra (graded coalgebra) we mean a non-negatively graded
vector space $A=\bop_{n=0}^\infty A_n$ ($C=\bop_{n=0}^\infty C_n$)
over~a~field~$\k$  such that $A_0=\k$ ($C_0=\k$) equipped with an
associative algebra (coalgebra) structure that respects the grading, i.~e.,
$A_i\cd A_j\sub A_{i+j}$ and $1\in A_0\,$
($\D(C_n)\sub\sum_{i+j=n}C_i\ot C_j$ and $\e(C_{>0})=0$).
 A graded algebra (coalgebra) structure induces an augmented algebra
(coalgebra) structure in an evident way.
 The homology coalgebra (cohomology algebra) of a graded algebra
(coalgebra) is equipped with a natural second grading, as it can be seen
from the explicit resolutions above:
  $$
   H_*(A)=\bop_{i\le j}H_{ij}(A)
   \quad \text{and} \quad
   H^*(C)=\bop_{i\le j}H^{ij}(C).
  $$

 In fact, all the results below in this section can be formulated in a more
general setting of a graded algebra in a (semisimple abelian, not
necessarily symmetric) tensor category, where the duality connects the
algebras in the opposite categories; however, we prefer to deal with vector
spaces here.

\subsection{Quadratic algebras and coalgebras} \

\begin{rem}{Definition 1.}
 A graded coalgebra $C$ is called {\it one-cogenerated\/} if the iterated
comultiplication maps $\D^{(n)}\:C_n\rarrow C_1^{\ot n}$ are injective,
or equivalently, all the maps $\D\:C_{i+j}\rarrow C_i\ot C_j$ are injective.
 A graded coalgebra is called {\it quadratic\/} if it is isomorphic to the
subcoalgebra of the form
 $$
  \lan V,R\ran=\bop_{n=0}^\infty\,\bcap_{i=1}^{n-1} V^{i-1}\ot R\ot V^{n-i-1}
 $$
of the tensor coalgebra $\N(V)=\bop_n V^{\ot n}$ for some vector space $V$
and a subspace $R\sub V^{\ot2}$.
 With a graded coalgebra $C$, one can associate in a natural way a
quadratic coalgebra $\q C$ and a morphism of graded coalgebras
$r_C\:C\rarrow \q C$ that is an isomorphism on $C_1$ and an epimorphism
on $C_2$.
 A graded algebra is called {\it quadratic\/} if it is isomorphic
to the quotient algebra $\{V,R\}=\T(V)/(R)$ of a tensor algebra
$\T(V)=\bop_n V^{\ot n}$ by the ideal generated by a subspace
$R\sub V^{\ot2}$.
 With a graded algebra $A$, one can associate
a quadratic algebra $\q A$ and a morphism of graded algebras
$r_A\:\q A\rarrow A$ that is an isomorphism on $A_1$ and a monomorphism
on $A_2$.
 \end{rem}

\begin{remsec}{Definition 2.}
 The quadratic algebra $A=\{V,R\}$ and the quadratic coalgebra
$C=\lan V,R\ran$ are called {\it dual\/} to each other; we denote this as
$C=A^!$ and $A=C^?$.
 Evidently, this defines an equivalence between the categories of
quadratic algebras and quadratic coalgebras.
\end{remsec}

\begin{thmsec}{Proposition 1}
 A graded coalgebra $C$ is one-cogenerated iff one has $H^{1,j}(C)=0$ for
$j>1$.
 A one-cogenerated coalgebra $C$ is quadratic iff $H^{2,j}(C)=0$ for $j>2$.
 Moreover, the morphism $r_C\:C\rarrow\q C$ is an isomorphism on
$C_{\le n}$ iff $H^{2,j}(C)=0$ for $2<j\le n$.
 The analogous statements are true for graded algebras.
\end{thmsec}

\pr{Proof}:
 It is evident from the explicit form of the cobar-complex that
$H^{1,>1}(C)=0$ for a one-cogenerated coalgebra $C$.
 Conversely, let $j>1$ be the minimal number for which
$\D^{(j)}\:C_j\rarrow C_1^{\ot j}$ is non-injective, then it is easy to see
that the map $\D\:C_j\rarrow \bop^{s+t=j}_{s,t\ge1} C_s\ot C_t$ is
non-injective also, hence $H^{1,j}(C)\ne0$.
 Now let $C$ be one-cogenerated, then the map $r_C$ is an embedding.
 Let $z\in C_+\ot C_+$ be a homogeneuos cocycle of degree $n$,
thus $z=\sum^{s+t=n}_{s,t\ge1} z_{st}$, where $z_{st}\in C_s\ot C_t$.
 The cocycle condition means that the images of $(\D\ot\id)(z_{u+v,w})$ and
$(\id\ot\D)(z_{u,v+w})$ in $C_u\ot C_v\ot C_w$ coincide for any
$u,v,w\ge1$, $\,u+v+w=n$.
 Since the maps $\D^{(k)}$ are injective, it is equivalent to say that
the elements $(\D^{(s)}\ot\D^{(t)})(z_{st})\in C_1^{\ot n}$ coincide for
all $s$ and $t$.
 We have got an element in $C_1^{\ot n}$; it is easy to see that it
represents an element of $\q C$ which belongs to the image of $r_C$ iff $z$
is a coboundary.
 At last, if $r_C$ is an isomorphism in degree $<n$, then any element of
$\q C$ corresponds to a cocycle $z$ in this way.
 \qed

\begin{thm}{Proposition 2}
 For any graded coalgebra $C$, the diagonal subalgebra
$\bop_i H^{i,i}(C)$ of the cohomology algebra $H^*(C)$ is
a quadratic algebra isomorphic to $(\q C)^?$.
 Analogously, for a graded algebra $A$, the diagonal quotient coalgebra
$\bop_i H_{i,i}(A)$ of the homology coalgebra $H_*(A)$ is
isomorphic to $(\q A)^!$.
\end{thm}

\pr{Proof} is immediate. \qed

\subsection{Koszul algebras and coalgebras}
 This definition is due to S.~Priddy~\cite{Pr};
see also~\cite{Lof,Bac,BF,BGS}.

\begin{rem}{Definition 3.}
 A graded algebra $A$ is called {\it Koszul\/} if $H_{ij}(A)=0$ unless
$i=j$.
 A graded coalgebra $C$ is called {\it Koszul\/} if $H^{ij}(C)=0$ unless
$i=j$.
 It follows from Proposition~1 that any Koszul algebra (coalgebra) is
quadratic.
 \end{rem}

 Now we are going to establish the criterion of Koszulity in the explicit
linear algebra terms due to J.~Backelin~\cite{Bac}.
 In particular, we will see that the dual algebra and coalgebra are Koszul
simultaneuosly.

\begin{rem}{Definition 4.}
 A collection of subspaces $X_1$, \ds, $X_{n-1}$ in a vector space
$W$ is called {\it distributive}, if there exists a (finite)
direct decomposition $W=\bop_{\o\in\O}W_\o$ such that each
subspace $X_k$ is the sum of a set of subspaces $W_\o$.
 Equivalently, the distributivity identity $(X+Y)\cap Z=X\cap Z+Y\cap Z$
should be satisfied for any triple of subspaces $X$, $Y$, $Z$ that can be
obtained from the subspaces $X_k$ using the operations of sum and
intersection.
 \end{rem}

\begin{thmsec}{Lemma}
 Let $X_1$, \ds, $X_{n-1}\sub W$ be a collection of linear subspaces;
assume that any its proper subcollection
$X_1$, \ds, $\widehat{X}_k$, \ds, $X_{n-1}$ is distributive.
 Then the following three conditions are equivalent:
  \begin{enumerate}
    \item[(a)] the following complex $B^*(W,X)$ is exact everywhere outside
               its left term:
            $$
             W\rarrow \bop_s W/X_s \rarrow \bop_{s<t}W/(X_s+X_t)
             \rarrow \ds \rarrow W/\textstyle\sum_k X_k \rarrow0;
            $$
    \item[(b)] the following complex $B_*(W,X)$ is exact everywhere outside
               its left term:
            $$
             W\larrow \bop_s X_s \larrow \bop_{s<t}X_s\cap X_t
             \larrow \ds \larrow \textstyle\bcap_k X_k \larrow0;
            $$
    \item[(c)] the collection $X_1$, \ds, $X_{n-1}$ is distributive. \qed
  \end{enumerate}
\end{thmsec}

\begin{thmsec}{Proposition 3}
 Let $V$ be a vector space and $R\sub V\ot V$ be a subspace, then the
following three conditions are equivalent:
  \begin{enumerate}
    \item[(a)] the quadratic algebra $A=\{V,R\}$ is Koszul;
    \item[(b)] the quadratic coalgebra $C=\lan V,R\ran$ is Koszul;
    \item[(c)] for any $n$, the collection of subspaces
               $V^{\ot k-1}\ot R\ot V^{n-k-1}\sub V^{\ot n}$,
               where $k=1$,~\ds,~$n-1$, is distributive.
  \end{enumerate}
\end{thmsec}

\pr{Proof}:
 Moreover, one has $H_{ij}(A)=0$ for $i<j\le n$ iff the collection of
subspaces in $V^{\ot n}$ is distributive, and the same for coalgebras.
 This follows immediately from Lemma by induction on $n$.  \qed

\Section{Cohomology of Nilpotent Coalgebras}

\subsection{Nilpotent coalgebras}
 Let $C$ be an augmented coalgebra with the augmentation map
$\c\:\k\rarrow C$.
 The {\it augmentation filtration\/} on an augmented coalgebra $C$ is an
increasing filtration $N$ defined by the formula
  $$
   N_nC=\bs c\in C \mid \D^{(n+1)}(c)\in C^{\ot n+1}_\c = \sum_{i=1}^{n+1}
C^{\ot i-1}\ot\c(\k)\ot C^{\ot n-i+1}\sub C^{\ot n+1} \es,
  $$
where $\D^{(m)}\:C\rarrow C^{\ot m}$ denotes the iterated comultiplication
map.
 In particular, we have $N_0C=\c(\k)$.

\begin{thm}{Proposition 4}
 The filtration $N$ respects the coalgebra structure on $C$, that is
  $$
   \D(N_nC)\sub \sum_{i+j=n} N_iC\ot N_jC.
  $$
 Furthermore, the associated graded coalgebra
$\gr_NC=\bop_{n=0}^\infty N_nC/N_{n-1}C$ is one-cogenerated.
 \end{thm}

\pr{Proof}:
 Let $\phi\:C\rarrow\k$ be a linear function annihilating $N_{k-1}C$,
where $0\le k\le n$; then it can be factorized as
$\phi=\psi\circ\D^{(k)}$,
where
$\psi\:C^{\ot k}\rarrow\k$
is a function annihilating
$C^{\ot k}_\c$.
 We have to show that
$(\phi\ot\nobreak\id)\D N_nC\sub N_{n-k}C$.
 Put for convenience
$\D^{(0)}=\e$ and $\D^{(1)}=\id$;
then one has
$(\D^{(k)}\ot\nobreak\D^{(n-k+1)})\circ\D=\D^{(n+1)}$,
hence
$(\phi\ot\nobreak\D^{(n-k+1)})\D N_nC=
(\psi\ot\nobreak\id^{\ot n-k+1})\D^{(n+1)}N_nC \sub C^{\ot n-k+1}_\c$,
so we are done.
 Since we have $\D^{(n)}(c)\notin C^{\ot n}_\c$ for $c\notin N_{n-1}C$, the
second assertion is immediate.
 \qed

\begin{rem}{Definition 5.}
 An augmented coalgebra $C$ is called {\it nilpotent\/} if the augmentation
filtration $N$ is full, that is $C=\bigcup_nN_nC$.
 \end{rem}

\begin{remsec}{Example:}
 Let $G$ be a pro-$l$-group and $C=\F_l(G)$ be the coalgebra of locally
constant functions on $G$ with respect to the convolution; in other words,
$C=\varinjlim\F_l(G/U)$, where the limit is taken over all open
normal subgroups $U$ of $G$ and the coalgebra $\F_l(G/U)=\F_l[G/U]^*$ is
the dual vector space to the group algebra of~$G/U$.
 Let $\c:\F_l\rarrow\F_l(G)$ be the augmentation map that takes a constant
from $\F_l$ to the corresponding constant function on $G$.
 Since the augmentation ideal of the group ring of a finite $l$-group
over $\F_l$ is nilpotent, it follows by passing to the inductive
limit that the augmented coalgebra $\F_l(G)$ is nilpotent also.
 It is easy to see that the category of $\F_l(G)$-comodules is equivalent
to the category of discrete $G$-modules over $\F_l$ (and the same is true
over $\Z$) for any pro-finite group $G$.
 \end{remsec}

\subsection{Main theorem}
 Now we are ready to prove the theorem mentioned in Introduction.

\begin{thm}{Theorem}
 Let $H=H^*(C)$ be the cohomology algebra of a nilpotent coalgebra $C$.
 Assume that
  \begin{enumerate}
    \item $H^2$ is generated by $H^1$;
    \item in the subalgebra generated by $H^1$ in $H$, there are no
          nontrivial relations of degree~3;
    \item the quadratic algebra $\q H$ defined by $H^1$ and $H^2$ is
          Koszul.
  \end{enumerate}
 Then the whole algebra $H$ is quadratic (and therefore, Koszul).
 In addition, there is an isomorphism $H^*(C)\simeq H^*(\gr_NC)$.
\end{thm}

\pr{Proof}:
 The filtration $N$ on a coalgebra $C$ induces a filtration on the
corresponding cobar-complex:
  $$
   N_nC^{+\ot i}=\bop_{j_1+\ds+j_i=n}N_{j_1}C^+\ot\ds\ot N_{j_i}C^+,
  $$
where $N_jC^+=N_jC/\c(\k)$, so that the filtration on $C^{+\ot i}$ starts
with $N_i$.
 Clearly, the associated graded complex coincides with the cobar-complex of
$\gr_NC$, thus we obtain a multiplicative spectral sequence
  $$
   E_1^{ij}=H^{ij}(\gr_NC) \implies H^i(C),
  $$
which converges since the filtration is an increasing one.
 More exactly, the differentials have the form
$d_r\:E_r^{i,j}\rarrow E_r^{i+1,j-r}$
and there is an induced increasing multiplicative filtration $N$ on
$H^*(C)$ such that $\gr_N^jH^i(C)=E_\infty^{i,j}$.

 In particular, we see that the subalgebra $\bop N_iH^i(C)$ in $H^*(C)$
is isomorphic to the quotient algebra of the diagonal cohomology
$\bop H^{i,i}(\gr_NC)$ by the images of the differentials.
 By  Propositions~4, the graded coalgebra $\gr_NC$ is one-cogenerated,
hence (by Proposition~1) we have $E_1^{1,j}=H^{1,j}(\gr_NC)=0$ for $j>1$,
which implies $H^1(C)=N_1H^1(C)\simeq H^{1,1}(\gr_NC)$ and
$N_2H^2(C)\simeq H^{2,2}(\gr_NC)$.
 Since (by Proposition~2) the diagonal cohomology algebra
$\bop{}H^{i,i}(\gr_NC)$ is quadratic, we conclude that it is isomorphic
to~$\q H^*(C)$.
 By Proposition~2 again, $\bop H^{i,i}(\gr_NC)$ is the dual quadratic
algebra to the coalgebra $\q\gr_NC$; since we suppose $\q H^*(C)$ is
Koszul, the dual coalgebra $\q\gr_NC$ is Koszul also (Proposition~3).
 On the other hand, we have assumed that there are no cubic relations in
the subalgebra generated by $H^1(C)$, hence all the differentials
$d_r\:E_r^{2,3+r}\rarrow E_r^{3,3}$ targeting in $H^{3,3}(\gr_NC)$ vanish.

 Now let us prove by induction that $H^{2,j}(\gr_NC)=0$ for $j>2$.
 Assume that this is true for $2<j\le n-1$; by Proposition~1, it follows
that the map $r_{\gr_NC}\:\gr_NC\rarrow\q\gr_NC$ is an isomorphism in
degree $\le n-1$.
 Therefore, the induced map on the cobar-complex is an isomorphism in
these degrees also, hence in particular $H^{3,j}(\gr_NC)=H^{3,j}(\q\gr_NC)$
for $j\le n-1$ (and even for $j\le n$).
 Since the coalgebra $\q\gr_NC$ is Koszul, it follows that
$E_1^{3,j}=H^{3,j}(\gr_NC)=0$ for $3<j\le n-1$ and the term
$E_1^{2,n}=\allowbreak H^{2,n}(\gr_NC)$ cannot die in the spectral sequence.
 But we have assumed that $H^2(C)$ is generated by $H^1(C)$, hence
$H^2(C)=N_2H^2(C)$ and $E_\infty^{2,n}=0$, so we are done.

 We have seen that the coalgebra $\gr_NC$ is quadratic and $\q\gr_NC$ is
Koszul, that is $\gr_NC$ is Koszul.
 It follows that $E_1^{i,j}=0$ for $i\ne j$, thus the spectral sequence
degenerates and $H^*(C)=H^*(\gr_NC)$.
 Therefore, $H^*(C)$ is Koszul also.
\qed

\begin{rem}{Remark:}
 This result is a formal analogue of some kind of
Poincare--Birkhoff--Witt theorem for filtrations on quadratic
algebras~\cite{PP}; in other words, it can be considered as reflecting
the deformation properties of Koszul algebras.
\end{rem}

 In the conclusion, recall the consequences we get for the Bloch--Kato
conjecture.
 Since the conditions (1) and (2) of our Theorem are known to be satisfied
for the coalgebra $C=\F_2(G_F)$ of any absolute Galois pro-2-group $G_F$
and the quadratic part $\q H^*(C)$ of the corresponding cohomology algebra
is exactly the Milnor K-theory algebra $\KM(F)\ot\F_2$, it suffices to
establish the Koszul property of this quadratic algebra in order to prove
the conjecture for $l=2$.
 The same would be true for the other $l$ if we know the norm
residue homomorphism for that $l$ to be injective in degree~$3$.

\clearpage

\bigskip

\appendix

\section*{Appendix: $\protect\KM(F)\ot\F_l$ is Koszul for
  all Primitive Fields $F$}
\setcounter{section}{1}

\medskip

\subsection{Commutative PBW-bases}
 Let $A$ be a commutative or skew-commutative one-generated algebra over a
field $\k$.
 Fix a basis $\X=\{x\}$ of the vector space $A_1$; suppose
$\X$ is equipped with a complete order $x'<x''$, i.~e., there are no
infinite decreasing sequences $x_1>x_2>x_3>\dsb\;$.
 Introduce the inverse lexicographical order on the set $\X^{(n)}$ of all
monomials in $\X$ of fixed degree~$n\/$:
$\prod_x x^{i_x}<\prod_xx^{j_x}$,
where $\sum_xi_x=\sum_xj_x=n$,
if there is $x_0\in\X$ such that $i_x=j_x$
for all $x<x_0$ and $i_{x_0}<j_{x_0}$.
 It easy to see that this is a complete order also.

\begin{thm}{Lemma}
 Let $W$ be a vector space spanned by a completely ordered set of its
vectors $w\in\W$.
 Then the set of all $w\in\W$ which cannot be expressed as a (finite)
linear combination of the smaller ones forms a basis of $W$. \qed
 \end{thm}

 Applying this statement to the space $W=A_n$ generated by the images of
monomials of degree $n$, we obtain a monomial basis $S_n\sub\X^{(n)}$
in $A_n$.
 It is clear that any monomial of degree $k$ which divides a monomial from
$S_n$ belongs to $S_k$.
 This basis of $A$ is called a {\it commutative\/} PBW-{\it basis\/}
iff the set $S_n$ coincides with the set of all monomials of degree
$n$ whose divisors of degree~$2$ belong to $S_2$.
 It is easy to see that any algebra admitting a commutative PBW-basis is
quadratic.

\begin{rem}{Remark:}
 In the case of a (skew-)commutative algebra, the notion of a commutative
PBW-basis is more general then the well-known Priddy's definition~\cite{Pr}
of a (non-commutative) PBW-basis.
 The following results are completely analogous to the non-commutative
case.
\end{rem}

\begin{thmsec}{Proposition}
 Assume that the PBW condition above is satisfied for the monomial basis
$S_3$ in the third degree component $A_3$ of a quadratic algebra $A$.
 Then the same is true for any degree~$n$.
 Furthermore, a quadratic algebra admitting a commutative PBW-basis is
Koszul.
\end{thmsec}

\pr{Proof} can be found in~\cite{PP}.
 The first statement is well-known~\cite{aBer,aBuc}.
 The second one was proved in~\cite{Kem} using the result of
R.~Fr\"oberg~\cite{aFr}.
 Some of these papers deal with the finite dimensional case,
but this is not essential here.
 \qed

\subsection{Finite fields}
 Since $\KM_2(F)=0$ for any finite field $F$ \cite{aMil1}, this ring is
evidently Koszul.

\subsection{The field of rational numbers}
 Let us show that the algebra $\KM(\Q)\ot\F_l$ is Koszul for any prime
number $l$.
 (This will not be applicable to our problem directly since $G_\Q$ is not
a pro-$l$-group for any $l$.)
 By the well-known results of J.~Milnor~\cite{aMil2} and H.~Bass and
J.~Tate~\cite{aBT}, one has
 $$
  \KM_n(\Q)\trarrow\bop_{\text{prime }p}
  \KM_{n-1}(\F_p)\;\,{\textstyle\bop}\;\,\KM_n(\R)\ot\F_2,
 $$
where $\d_p\:\KM_n(\Q)\rarrow\KM_{n-1}(\F_p)$ is a boundary homomorphism.
 Explicitly, we have
  \begin{align*}
    &\KM_2(\Q)\ot\F_l\trarrow\bop_{\text{prime }p>2}\F_p^*/(F_p^*)^l
      \;\,{\textstyle\bop}\;\,
      \cases\Z/2\Z,\quad\text{for \ }l=2\\
      0,\quad\text{for \ }l>2\endcases \\
    &\KM_n(\Q)\ot\F_l\trarrow
      \cases\Z/2\Z,\quad\text{for \ }l=2\\
      0,\quad\text{for \ }l>2\endcases
  \end{align*}
for $n\ge3$, where the multiplication is given by the formulas
  \begin{align*}
    &\{a,b\}_p=(-1)^{\nu_p(a)\nu_p(b)}\,
      \frac{a^{\nu_p(b)}}{b^{\nu_p(a)}}\mod p \\
    &\{a_1,\ds,a_n\}_\infty=
      \cases1\mod2,\quad\text{if \ }a_1,\ds,a_n<0\\
      0\mod2,\quad\text{otherwise.}\endcases
  \end{align*}
 Let us construct PBW-bases for these algebras in an explicit way.
 Consider two cases separately.

(a) $l=2$.
 The vector space $\KM_1(\Q)\ot\F_2$ admits a basis consisting of the
symbols $\{p\}$, for all prime numbers $p$, and $\{-1\}$.
 Let $R$ be the set of all odd primes $r$ for which $2$ is a square modulo
$r$ and $Q$ be the set of all other odd primes $q$.
 Fix arbitrary complete orders on the sets $Q$ and $R$ and extend them
to an order on the whole basis by the rule $\{2\}<\{q\}<\{r\}<\{-1\}$
for any $q\in Q$ and $r\in R$.
 For each $r\in R$, let $q(r)$ denote the minimal prime $q\in Q$ with
respect to the order which is not a square modulo $r$.
 Since $\{2,2\}=\{2,r\}=0$ for any $r\in R$, it is easy to verify
that the monomials $\{2,q\}$, $\,\{q(r),r\}$, and $\{-1,\ds,-1\}$ form
the basis of $\KM(\Q)\ot\F_2$ corresponding to the order chosen.
 The PBW condition is clearly satisfied.

(b) $l>2$.
 The vector space $\KM_1(\Q)\ot\F_l$ admits a basis consisting of the
symbols $\{p\}$ for all prime numbers $p$.
 Let $R$ be the set of all primes $r$ for which $r-1$ is divided by $l$
and $Q$ be the set of all other primes $q$.
 Choose a complete order on our basis such that $\{q\}<\{r\}$ for any
$q\in Q$ and $r\in R$.
 For each $r\in R$, let $q(r)$ denote the minimal $q\in Q$ which is not
a $l$-th power modulo $r$.
 Since $\{q,q'\}=0$ and $\{q,r\}_q=0$ for any $q$, $q'\in Q$ and
$r\in R$, it is easy to see that the monomials $\{q(r),r\}$ form a
commutative PBW-basis of $\KM(\Q)\ot\F_l$.


\clearpage

\bigskip

\begin{thebibliography}{9}

\smallskip

\bibitem{Bac}
 J.~Backelin.
  A distributiveness property of augmented algebras and some
related homological results.
  Ph.~D. Thesis, Stockholm, 1981.

\bibitem{BF}
 J.~Backelin, R.~Fr\"{o}berg.
  Koszul algebras, Veronese subrings and and rings with linear
resolutions.
{\it Rev. Roumaine Math. Pures Appl.} {\bf 30}, \#2, p.~85--97, 1985.

\bibitem{aBT}
 H.~Bass, J.~Tate.
   The Milnor ring of a global field.
{\it Lecture Notes in Math.} {\bf 342}, p.~349--446, 1973.

\bibitem{BGS}
 A.~A.~Beilinson, V.~A.~Ginzburg, V.~V.~Schechtman.
  Koszul duality.
{\it Journ. Geom. Phys.} {\bf 5}, \#3, p.~317--350.

\bibitem{aBer}
 G.~Bergman.
   The diamond lemma for rings theory.
{\it Advances in Math.} {\bf 29}, \#2, p.~178--218, 1983.

\bibitem{aBuc}
 B.~Buchberger.
  Grobner bases: An algorithmic method in polinomial ideal
theory.
{\it CAMP--Bull.} {\bf 290}, \#83, 1978.

\bibitem{aFr}
 R.~Fr\"{o}berg.
   Determination of a class of Poincare series.
{\it Math. Scand.} {\bf 37}, p.~29--39, 1975.

\bibitem{Lof}
 C.~L\"{o}fwall.
   On the subagebra generated by the one-dimensional elements
in the Yoneda Ext-algebra.
{\it Lecture Notes in Math.} {\bf 1183}, p.~291--338, 1986.

\bibitem{Kem}
 G.~R.~Kempf. Some wonderful rings in algebraic geometry.
{\it Journ. Algebra\/} {\bf 134}, p.~222--224, 1990.

\bibitem{MS1}
 A.~S.~Merkurjev, A.~A.~Suslin.
    $K$-cohomology of Severi--Brauer varieties and the norm residue
homomorphism.
{\it Math. USSR Izvestiya\/} {\bf 21}, \#2, p.~307--340, 1983.

\bibitem{MS2}
 A.~S.~Merkurjev, A.~A.~Suslin.
    The norm residue homomorphism of degree three.
{\it Math. USSR Izvestiya\/} {\bf 36}, \#2, p.~349--367, 1991.

\bibitem{aMil1}
 J.~Milnor.
   Introduction to algebraic $K$-theory.
{\it Ann. Math. Studies}, Princeton, 1971.

\bibitem{aMil2}
 J.~Milnor.
   Algebraic $K$-theory and quadratic forms.
{\it Invent. Math.} {\bf 9}, p.~318--344, 1970.

\bibitem{PP}
 A.~Polishchuk, L.~Positselski.
   Quadratic algebras.
In preparation.

\bibitem{Pr}
 S.~Priddy.
   Koszul resolutions.
{\it Trans. AMS\/} {\bf 152}, p.~39--60, 1970.

\bibitem{Ros}
 M.~Rost.
   Hilbert theorem~90 for $K_3$ for degree two extensions.
Preprint, 1986.

\end{thebibliography}
\end{document}